# On the Jensen functional and superterzaticity

**MITROI-SYMEONIDIS Flavia-Corina**

University of South-East Europe LUMINA, Faculty of Engineering Sciences, Şos. Colentina 64B, RO-021187, Bucharest, Romania

**Email address:**
fcmitroi@yahoo.com

**Abstract:** In this note we describe some results concerning upper and lower bounds for the Jensen functional. We use several known and new results to shed light on the concepts of superterzatic functions.

## 1. Introduction

The aim of this paper is to discuss new results concerning the Jensen functional in the framework of superterzatic functions.

For the reader's convenience, before going into details, we quote here some relevant results regarding the superterzaticity and the Jensen functional.

**Definition 1.** We consider a real valued function *f* defined on an interval *I*, $x_1, x_2, ..., x_n \in I$ and $p_1, p_2, ..., p_n \in (0,1)$ with $\sum_{i=1}^{n} p_i = 1$. The *Jensen functional* is defined by

$$\mathrm{J}(f,\mathbf{p},\mathbf{x}) = \sum_{i=1}^{n} p_i f(x_i) - f\left(\sum_{i=1}^{n} p_i x_i\right)$$

(see [2]).

**Definition 2**. A function *f* defined on an interval $I=[0,a)$ is superterzatic, if for each $\bar{x} \in I$ there exists a real number $\mathrm{C}(\bar{x})$ such that

$$\mathrm{J}(f,\mathbf{p},\mathbf{x}) \geq \sum_{i=1}^{n} p_i x_i \left[(x_i-\bar{x})\mathrm{C}(\bar{x}) + \frac{f(|x_i-\bar{x}|)}{|x_i-\bar{x}|}\right]$$

$$= \sum_{i=1}^{n} p_i (x_i-\bar{x})^2 \mathrm{C}(\bar{x}) + \sum_{i=1}^{n} p_i x_i \frac{f(|x_i-\bar{x}|)}{|x_i-\bar{x}|}$$

for all $x_i \in I$, $i=1,...,n$, and $p_i \geq 0$, $i=1,...,n$, with $\sum_{i=1}^{n} p_i = 1$, such that $\bar{x}=\sum_{i=1}^{n} p_i x_i$.

We use the convention f(0)/0=0.

This definition was mentioned by S. Abramovich in her talk given at the Conference on Inequalities and Applications 10, [1].

The set of superterzatic functions is closed under addition and positive scalar multiplication.

**Example 3**. [1] Let $f(x) = x^p$, $p \geq 3$. This function is superterzatic with $\mathrm{C}(\bar{x}) = p\bar{x}^{p-1}$. The function $g(x)=x^3 \log(1/x)$, $g(0)=0$ is also subterzatic (i.e. the inequality in Definition 2 holds with reversed sign).

In what follows we shall be also interested in a more general Jensen functional and its behavior in the context of superterzaticity.

## 2. Main results

We introduce in a natural way a more general functional.

**Definition 4**. Assume that we have a real valued function *f* defined on an interval *I*, the real numbers $p_{ij}$, $i=1,...,k$ and $j=1,...,n_i$ such that $p_{ij}>0$, $\sum_{j=1}^{n_i} p_{ij}=1$ for all i=1,...,*k* (we put $\mathbf{p}_i=(p_{i1},\ldots,p_{in_i})$), $\mathbf{x}_i=(x_{i1},\ldots,x_{in_i}) \in I^{n_i}$ for all i=1,...,*k* and $\mathbf{q}=(q_1,\ldots,q_k)$, $q_k>0$ such that $\sum_{i=1}^{k} q_i=1$. We define the generalized Jensen functional by

$$\mathrm{J}_k(f,\mathbf{p}_1,\ldots,\mathbf{p}_k,\mathbf{q},\mathbf{x}_1,\ldots,\mathbf{x}_k) = \sum_{j_1,\ldots,j_k=1}^{n_1,\ldots,n_k} p_{1j_1}\ldots p_{kj_k} f\left(\sum_{i=1}^{k} q_i x_{ij_i}\right) - f\left(\sum_{i=1}^{k} q_i \sum_{j=1}^{n_i} p_{ij} x_{ij}\right).$$

We notice that for *k*=1 this definition reduces to Definition 1.

For more results concerning Jensen's functional the reader is referred to the papers [4,5].

Before announcing the main result, let us give the following lemma that describes the behaviour of the functional under the superterzaticity condition:

**Lemma 5**. Let $f$, $\mathbf{p}_i, \mathbf{q}, \mathbf{x}_i$ be as in Definition 4. If *f* is superterzatic then we have

$$J_k(f, \mathbf{p}_1, \dots, \mathbf{p}_k, \mathbf{q}, \mathbf{x}_1, \dots, \mathbf{x}_k)$$

$$\geq \sum_{j_1,\dots,j_k=1}^{n_1,\dots,n_k} p_{1j_1} \dots p_{kj_k} \sum_{i=1}^{k} q_i x_{ij_i} \times$$

$$\times \left[ \left( \sum_{i=1}^{k} q_i x_{ij_i} - \bar{x} \right) \mathrm{C}(\bar{x}) + \frac{f\left( \left| \sum_{i=1}^{k} q_i x_{ij_i} - \bar{x} \right| \right)}{\left| \sum_{i=1}^{k} q_i x_{ij_i} - \bar{x} \right|} \right],$$

where $\bar{x} = \sum_{i=1}^{k} q_i \ \sum_{j=1}^{n_i} p_{ij} x_{ij}$.

**Proof.** Since

$$\sum_{j_1,\dots,j_k=1}^{n_1,\dots,n_k} p_{1j_1} \dots p_{kj_k} = 1$$

and

$$\sum_{j_1,\dots,j_k=1}^{n_1,\dots,n_k} p_{1j_1} \dots p_{kj_k} \sum_{i=1}^{k} q_i x_{ij_i} = \sum_{i=1}^{k} q_i \sum_{j=1}^{n_i} p_{ij} x_{ij}$$

we use Definition 2 and the conclusion follows.□

**Theorem 6**. Let *f*, $\mathbf{p}_i, \mathbf{q}, \mathbf{x}_i$ be as in Definition 4 and let the positive real numbers $r_{ij}$, $i$=1,...,$k$, and $j$=1,...,$n_i$, be such that $r_{ij}$>0, $\sum_{j=1}^{n_i} r_{ij}$=1 for all i=1,...,$k$. We put

$$\mathbf{p}_i = (p_{i1}, \dots, p_{in_i}) \text{ for all } i=1,\dots,k,$$

$$m = \min_{\substack{1 \le j_1 \le n_1 \\ \dots \\ 1 \le j_k \le n_k}} \left\{ \frac{p_{1j_1} \dots p_{kj_k}}{r_{1j_1} \dots r_{kj_k}} \right\},$$

$$M = \max_{\substack{1 \le j_1 \le n_1 \\ \dots \\ 1 \le j_k \le n_k}} \left\{ \frac{p_{1j_1} \dots p_{kj_k}}{r_{1j_1} \dots r_{kj_k}} \right\}.$$

If *f* is superterzatic then

$$J_k(f, \mathbf{p}_1, \dots, \mathbf{p}_k, \mathbf{q}, \mathbf{x}_1, \dots, \mathbf{x}_k) - mJ_k(f, \mathbf{r}_1, \dots, \mathbf{r}_k, \mathbf{q}, \mathbf{x}_1, \dots, \mathbf{x}_k) \geq$$

$$\sum_{j_1,\dots,j_k=1}^{n_1,\dots,n_k} (p_{1j_1} \dots p_{kj_k} - m r_{1j_1} \dots r_{kj_k}) \sum_{i=1}^{k} q_i x_{ij_i} \times$$

$$\times \left[ \left( \sum_{i=1}^{k} q_i x_{ij_i} - \bar{x} \right) \mathrm{C}(\bar{x}) + \frac{f\left( \left| \sum_{i=1}^{k} q_i x_{ij_i} - \bar{x} \right| \right)}{\left| \sum_{i=1}^{k} q_i x_{ij_i} - \bar{x} \right|} \right]$$

$$+ m \sum_{i=1}^{k} q_i \sum_{j=1}^{n_i} r_{ij} x_{ij} \times$$

$$\times \left[ \sum_{i=1}^{k} q_i \sum_{j=1}^{n_i} (r_{ij} - p_{ij}) x_{ij} \mathrm{C}(\bar{x}) + \frac{f\left( \left| \sum_{i=1}^{k} q_i \sum_{j=1}^{n_i} (r_{ij} - p_{ij}) x_{ij} \right| \right)}{\left| \sum_{i=1}^{k} q_i \sum_{j=1}^{n_i} (r_{ij} - p_{ij}) x_{ij} \right|} \right]$$

and

$$MJ_k(f, \mathbf{r}_1, \dots, \mathbf{r}_k, \mathbf{q}, \mathbf{x}_1, \dots, \mathbf{x}_k) - J_k(f, \mathbf{p}_1, \dots, \mathbf{p}_k, \mathbf{q}, \mathbf{x}_1, \dots, \mathbf{x}_k) \geq$$

$$\sum_{j_1,\dots,j_k=1}^{n_1,\dots,n_k} \left( M r_{1j_1} \dots r_{kj_k} - p_{1j_1} \dots p_{kj_k} \right) \sum_{i=1}^{k} q_i x_{ij_i} \times$$

$$\left[ \left( \sum_{i=1}^{k} q_i x_{ij_i} - \breve{x} \right) \mathrm{C}(\breve{x}) + \frac{f\left( \left| \sum_{i=1}^{k} q_i x_{ij_i} - \breve{x} \right| \right)}{\left| \sum_{i=1}^{k} q_i x_{ij_i} - \breve{x} \right|} \right]$$

$$+ \sum_{i=1}^{k} q_i \sum_{j=1}^{n_i} p_{ij} x_{ij} \times$$

$$\left[ \sum_{i=1}^{k} q_i \sum_{j=1}^{n_i} (p_{ij} - r_{ij}) x_{ij} \mathrm{C}(\breve{x}) + \frac{f\left( \left| \sum_{i=1}^{k} q_i \sum_{j=1}^{n_i} (p_{ij} - r_{ij}) x_{ij} \right| \right)}{\left| \sum_{i=1}^{k} q_i \sum_{j=1}^{n_i} (p_{ij} - r_{ij}) x_{ij} \right|} \right],$$

where $\bar{x} = \sum_{i=1}^{k} q_i \ \sum_{j=1}^{n_i} p_{ij} x_{ij}$ and $\breve{x} = \sum_{i=1}^{k} q_i \ \sum_{j=1}^{n_i} r_{ij} x_{ij}$.

**Proof.** *The first inequality.* Obviously

$$J_k(f, \mathbf{p}_1, \dots, \mathbf{p}_k, \mathbf{q}, \mathbf{x}_1, \dots, \mathbf{x}_k) - mJ_k(f, \mathbf{r}_1, \dots, \mathbf{r}_k, \mathbf{q}, \mathbf{x}_1, \dots, \mathbf{x}_k) =$$

$$\sum_{j_1,\dots,j_k=1}^{n_1,\dots,n_k} (p_{1j_1} \dots p_{kj_k} - m r_{1j_1} \dots r_{kj_k}) f\left( \sum_{i=1}^{k} q_i x_{ij_i} \right)$$

$$+ mf\left( \sum_{i=1}^{k} q_i \sum_{j=1}^{n_i} r_{ij} x_{ij} \right)$$

$$- f\left( \sum_{i=1}^{k} q_i \sum_{j=1}^{n_i} p_{ij} x_{ij} \right).$$

Since

$$\sum_{i=1}^{k} q_i \sum_{j=1}^{n_i} p_{ij} x_{ij}$$

$$= \sum_{j_1,\dots,j_k=1}^{n_1,\dots,n_k} (p_{1j_1} \dots p_{kj_k} - m r_{1j_1} \dots r_{kj_k}) \sum_{i=1}^{k} q_i x_{ij_i}$$

$$+ m \sum_{i=1}^{k} q_i \sum_{j=1}^{n_i} r_{ij} x_{ij},$$

the conclusion follows by Lemma 5.

*The second inequality.* One has

$$MJ_k(f, \mathbf{r}_1, \dots, \mathbf{r}_k, \mathbf{q}, \mathbf{x}_1, \dots, \mathbf{x}_k) - J_k(f, \mathbf{p}_1, \dots, \mathbf{p}_k, \mathbf{q}, \mathbf{x}_1, \dots, \mathbf{x}_k) =$$

$$\sum_{j_1,\dots,j_k=1}^{n_1,\dots,n_k} (M r_{1j_1} \dots r_{kj_k} - p_{1j_1} \dots p_{kj_k}) f\left( \sum_{i=1}^{k} q_i x_{ij_i} \right)$$

$$+ f\left( \sum_{i=1}^{k} q_i \sum_{j=1}^{n_i} p_{ij} x_{ij} \right)$$

$$- Mf\left( \sum_{i=1}^{k} q_i \sum_{j=1}^{n_i} r_{ij} x_{ij} \right).$$

Since

$$\breve{x} = \sum_{j_1,\dots,j_k=1}^{n_1,\dots,n_k} \frac{M r_{1j_1} \dots r_{kj_k} - p_{1j_1} \dots p_{kj_k}}{M} \sum_{i=1}^{k} q_i x_{ij_i}$$

$$+ \frac{1}{M} \sum_{i=1}^{k} q_i \sum_{j=1}^{n_i} p_{ij} x_{ij}$$

and

$$\sum_{j_1,\dots,j_k=1}^{n_1,\dots,n_k} \frac{M r_{1j_1} \dots r_{kj_k} - p_{1j_1} \dots p_{kj_k}}{M} + \frac{1}{M} = 1,$$

we may apply again Lemma 5 and the conclusion follows.□

The particular case $\mathbf{p}_1$=...=$\mathbf{p}_k$=$\mathbf{p}$ and $\mathbf{x}_1$=...=$\mathbf{x}_k$=$\mathbf{x}$ is of interest.

**Corollary 7**. We consider $\mathbf{x}=(x_1, x_2, ..., x_n)\in I^n$, $\mathbf{p}=(p_1, p_2, ..., p_n)\in(0,1)^n$ with $\sum_{i=1}^n p_i = 1$ and $\mathbf{r}=(r_1, r_2, ..., r_n)\in(0,1)^n$ with $\sum_{i=1}^n r_i = 1$. We put

$$m=\min_{1\le i_1,\dots,i_k\le n}\left\{\frac{p_{i_1}\dots p_{i_k}}{r_{i_1}\dots r_{i_k}}\right\},$$

$$M=\max_{1\le i_1,\dots,i_k\le n}\left\{\frac{p_{i_1}\dots p_{i_k}}{r_{i_1}\dots r_{i_k}}\right\}.$$

We define $\mathrm{J}_k(f,\mathbf{p},\mathbf{q},\mathbf{x})$

$$= \sum_{i_1,\dots,i_k=1}^{n} p_{i_1}\dots p_{i_k} f\left(\sum_{j=1}^{k} q_j x_{i_j}\right) - f\left(\sum_{i=1}^{n} p_i x_i\right).$$

If $f$ is superterzatic then

$\mathrm{J}_k(f,\mathbf{p},\mathbf{q},\mathbf{x})-m\mathrm{J}_k(f,\mathbf{r},\mathbf{q},\mathbf{x})\ge$

$$\sum_{i_1,\dots,i_k=1}^{n}(p_{i_1}\dots p_{i_k} - m r_{i_1}\dots r_{i_k})\sum_{j=1}^{k} q_j x_{i_j}$$

$$\times\left[\left(\sum_{j=1}^{k} q_j x_{i_j} - \bar{x}\right)\mathrm{C}(\bar{x}) + \frac{f\left(\left|\sum_{j=1}^{k} q_j x_{i_j} - \bar{x}\right|\right)}{\left|\sum_{j=1}^{k} q_j x_{i_j} - \bar{x}\right|}\right]$$

$$+m\sum_{i=1}^{n} r_i x_i$$

$$\times\left[\sum_{i=1}^{n}(r_i - p_i)x_i\,\mathrm{C}(\bar{x}) + \frac{f(|\sum_{i=1}^{n}(r_i - p_i)x_i|)}{|\sum_{i=1}^{n}(r_i - p_i)x_i|}\right]$$

and

$M\mathrm{J}_k(f,\mathbf{r},\mathbf{q},\mathbf{x})-\mathrm{J}_k(f,\mathbf{p},\mathbf{q},\mathbf{x})\ge$

$$\sum_{i_1,\dots,i_k=1}^{n}(M r_{i_1}\dots r_{i_k} - p_{i_1}\dots p_{i_k})\sum_{j=1}^{k} q_j x_{i_j}\times$$

$$\left[\left(\sum_{j=1}^{k} q_j x_{i_j} - \breve{x}\right)\mathrm{C}(\breve{x}) + \frac{f\left(\left|\sum_{j=1}^{k} q_j x_{i_j} - \breve{x}\right|\right)}{\left|\sum_{j=1}^{k} q_j x_{i_j} - \breve{x}\right|}\right]$$

$$+\sum_{i=1}^{n} p_i x_i\times$$

$$\left[\sum_{i=1}^{n}(p_i - r_i)x_i\,\mathrm{C}(\breve{x}) + \frac{f(|\sum_{i=1}^{n}(p_i - r_i)x_i|)}{|\sum_{i=1}^{n}(p_i - r_i)x_i|}\right].$$

where $\bar{x} = \sum_{i=1}^n p_i x_i$ and $\breve{x} = \sum_{i=1}^n r_i x_i$.

The case $k$=1 reduces nicely.

**Corollary 8**. We consider $x_1, x_2, ..., x_n\in I$, $p_1, p_2, ..., p_n\in(0,1)$ with $\sum_{i=1}^n p_i = 1$ and $r_1, r_2, ..., r_n\in(0,1)$ with $\sum_{i=1}^n r_i = 1$. We put

$$m=\min_{1\le i\le n}\left\{\frac{p_i}{r_i}\right\}, M=\max_{1\le i\le n}\left\{\frac{p_i}{r_i}\right\}.$$

If $f$ is a superterzatic function, then we have

$$\mathrm{J}(f,\mathbf{p},\mathbf{x}) - m\mathrm{J}(f,\mathbf{r},\mathbf{x}) \ge$$

$$\sum_{i=1}^{n}(p_i - m r_i)x_i\left[(x_i - \bar{x})\mathrm{C}(\bar{x}) + \frac{f(|x_i - \bar{x}|)}{|x_i - \bar{x}|}\right]$$

$$+m\sum_{i=1}^{n} r_i x_i\left[\sum_{i=1}^{n}(r_i - p_i)x_i\mathrm{C}(\bar{x}) + \frac{f(|\sum_{i=1}^{n}(r_i - p_i)x_i|)}{|\sum_{i=1}^{n}(r_i - p_i)x_i|}\right]$$

and

$$M\mathrm{J}(f,\mathbf{r},\mathbf{x}) - \mathrm{J}(f,\mathbf{p},\mathbf{x}) \ge$$

$$\sum_{i=1}^{n}(M r_i - p_i)x_i\left[(x_i - \breve{x})\mathrm{C}(\breve{x}) + \frac{f(|x_i - \breve{x}|)}{|x_i - \breve{x}|}\right]$$

$$+\sum_{i=1}^{n} p_i x_i\left[\sum_{i=1}^{n}(p_i - r_i)x_i\mathrm{C}(\breve{x}) + \frac{f(|\sum_{i=1}^{n}(p_i - r_i)x_i|)}{|\sum_{i=1}^{n}(p_i - r_i)x_i|}\right].$$